\documentclass{amsart}

\usepackage[backref=page]{hyperref}
\hypersetup{
  colorlinks=true, %set true if you want colored links
   linkcolor=red,  %choose some color if you want links to stand out
}

\newtheorem{theorem}{Theorem}[section]

\theoremstyle{definition}

\theoremstyle{remark}

\newtheorem{problem}[theorem]{Problem}
\numberwithin{equation}{section}

\begin{document}
\vspace{-.5in}
\title{
On the central role of the  scale invariant Poisson processes
on $(0,\infty) \ \ \ \ $
}
\title{
On the central role of scale invariant
Poisson~processes~on~$(0,\infty) \ \ \ \ $
}

% author one information
\author{Richard Arratia}

\curraddr{University of Southern California, Department of Mathematics, DRB 155, Los Angeles CA 90089-1113}
\email{ rarratia@math.usc.edu }
\thanks{ Work supported in part by NSF grant DMS
96-26412.}

\date{September 15, 1997; revision of January 16, 1998}

\begin{abstract}

The scale invariant Poisson processes on $(0,\infty)$ play a central
but mildly disguised
role in number theory, combinatorics, and genetics.  They  give
the continuous limits
which underly and unify diverse discrete structures, including
the prime factorization of a uniformly chosen integer,
the factorization of polynomials over finite fields,
the decomposition into cycles of random permutations, the decomposition
into components of random mappings, and the
Ewens sampling formula.  They deserve attention as one of
the fundamental and central objects of probability theory.

\end{abstract}

\def\beq{\begin{equation}}
\newcommand{\eq}{\begin{equation}}
\newcommand{\en}{\end{equation}}
\def\beqa{\begin{eqnarray}}
\def\ena{\end{eqnarray}}
\def\beqano{\begin{eqnarray*}}
\def\enano{\end{eqnarray*}}

%springer bbb characters E, P, 1
\def\bbbe{{\rm I\!E}}
\def\bbbp{{\rm I\!P}}

\def\bbbz{{\mathchoice {\hbox{$\sf\textstyle Z\kern-0.4em Z$}}
{\hbox{$\sf\textstyle Z\kern-0.4em Z$}}
{\hbox{$\sf\scriptstyle Z\kern-0.3em Z$}}
{\hbox{$\sf\scriptscriptstyle Z\kern-0.2em Z$}}}}

\def\bbbr{{\rm I\!R}}
\def\bbbn{{\rm I\!N}}

\def\bbbone{{\mathchoice {\rm 1\mskip-4mu l} {\rm 1\mskip-4mu l}
{\rm 1\mskip-4.5mu l} {\rm 1\mskip-5mu l}}}

\def\p{\bbbp}   %\def\p{{\bf P}}
\def\e{\bbbe}   %\def\e{{\bf E}}
\def\bone{\bbbone} %\def\bone{\mbox{\bf 1}}

\def\qed{\hfill \mbox{\rule{0.5em}{0.5em}}}
\def\ra{\rightarrow}
\def\Ra{{\Rightarrow}}
\def\no{\vskip 1pc \noindent}

\def\proof{\noindent{\bf Proof.\ \ }}
\def\indist{\mbox{$ \ =_d \ $}}
\def\inprob{\mbox{$\rightarrow_P\ $}}
\def\todist{\Rightarrow}
\def\bothsides{\stackrel{\mbox{\scriptsize $\smile$}}{\mbox{\scriptsize
$\frown$}}}

\def\qed{\mbox{\rule{0.5em}{0.5em}}}

\def\BP{\bbbp}
\def\BZ{\bbbz}
\def\BR{\bbbr}
\def\BN{\bbbn}
\def\V{{\mathcal V}}
\def\X{{\mathcal X}}
\def\L{{\mathcal L}}
\def\Y{{\mathcal Y}}
\def\D{{\mathcal D}}
\def\P{{\mathcal P}}
\def\L{{\mathcal L}}

\def\vec#1{\ifmmode
\mathchoice{\mbox{\boldmath$\displaystyle\bf#1$}}
{\mbox{\boldmath$\textstyle\bf#1$}}
{\mbox{\boldmath$\scriptstyle\bf#1$}}
{\mbox{\boldmath$\scriptscriptstyle\bf#1$}}\else
{\mbox{\boldmath$\bf#1$}}\fi}

\def\ba{\vec{a}}
\def\bb{\mbox{\boldmath $b$}}
\def\be{\mbox{\boldmath $e$}}
\def\bu{\mbox{\boldmath $u$}}
\def\bC{\mbox{\boldmath $C$}}
\def\bD{\mbox{\boldmath $D$}}
\def\bZ{\mbox{\boldmath $Z$}}

\def\th{\theta}

\def\note#1{{\medskip {\small \it #1} \medskip}}

\maketitle

\vspace{-.3in}
%\tableofcontents
\setcounter{section}{-1}
\section{Introduction}
The  scale invariant Poisson processes on
$(0,\infty)$ have almost always been overlooked, even though
these processes play a fundamental role in
combinatorics, number theory, and genetics.  Quantities and
relations which are most easily explained in terms of these Poisson
processes have been studied directly without mentioning the Poisson
connection.  Examples in number theory (\ref{dick buch})
include Dickman's function
$\rho$ from \cite{dickman} in 1930, Buchstab's function $\omega$
from \cite{buch} in 1937, and the ``convolution powers of
the Dickman function'' from \cite{db vl,hensley}.  The examples
from combinatorics and genetics involve the  Poisson-Dirichlet
process and the GEM process from the 1970's, which
  are usually discussed in terms
of products of independent, Beta-distributed random variables, or
constructed by conditioning a Poisson process which is not scale invariant.
We want to reveal the presence of the
scale invariant Poisson processes,
as well as their intrinsic beauty and simplicity,
for an audience including probabilists,
number theorists, geneticists, and others, without assuming background
such as knowledge of point processes.

The sections correspond, more or less, to the transparancies from the lecture, except for
the last two sections, which match discussions from the problem sessions, and have been substantially expanded.

\section{Scale invariant Poisson processes on $(0,\infty)$}

The only reference to ``scale invariant Poisson processes''
per se that we have found is Daley and Vere Jones \cite{point processes},
page 325, where the scale invariant Poisson processes on
$\BR$ are fully classified. We requested help from the audience in compiling a
bibliography on scale invariant Poisson processes; C. Newman supplied
a reference to \cite{chuck}.  Pitman and Yor \cite{pityor96} is a general study
of spacings in scale invariant random closed subsets of
 $(0,\infty)$, including examples
like the zeroes of Brownian motion, and the scale invariant Poisson processes.

The first subtle thing to realize about scale invariant
Poisson process{\em es} is that the plural is
essential! As with translation invariant Poisson processes,
there is a nonnegative parameter $\theta$ which appears as a
scalar factor for the intensity measure, but unlike the translation
invariant case,  varying $\theta$
leads to major qualitative changes.  In particular, for certain
natural aspects of qualitative behavior, there are
phase transitions at $\theta=1$ (see (\ref{super}),)
and at $\theta = 1 / \log 2$ (see section 22.)

Our notation for scaling is $cX := \{ cx: \ x \in X \}$,
defined for any $X \subset R^d$ and $c>0$.
For a random set $\X$, scale invariance is the
property that for all $c>0$, $c\X \indist \X$, i.e.~ all
scalings of the set have the same distribution.

In case $\X \subset (0,\infty)$, one can take logarithms:
let $\L =\log(\X) := \{ \log(x):\  x~\in~\X \}$.  In this
situation, clearly, scale invariance for $\X$ is equivalent to
translation invariance for $\L$.

\section{Intensity}
For us, a point process is a random set or multiset
 $\X$ such that for
sets $I$ in an appropriate class, $\X \cap I$ is
a finite set or multiset, with finite mean size; the points of
$\X$ are {\em not} labelled.
The intensity measure $\mu$ is then defined
by $\mu(I) = \e |\X \cap I|$, the expected number of points
falling in $I$.  For instance, if $\X$ is a translation invariant
point process on $\BR$, then its intensity is a translation
invariant, nonnegative measure, and hence has the form
$\theta$ times Lebesgue measure --- we will say that the
intensity is $\theta \ dx$ on $\BR$.
   The interesting intensity measures for our purposes are, each
with parameter $\theta > 0$:
$$
\theta \ dx \mbox{ on } (-\infty,\infty) \ \ \ \  \mbox{ translation invariant},
$$
$$
(\theta / x)  \ dx \mbox{ on } (0,\infty) \ \ \ \  \mbox{ scale invariant},
$$
$$
(\theta e^{-x} /x) \ dx \mbox{ on } (0,\infty) \ \ \ \    \mbox{``Gamma'' or ``Moran''
subordinator}.
$$
One can easily verify that
$(\theta/u) \ du $ gives the intensity of scale invariant point
processes $\X$ on $(0,\infty)$, starting from  the translation
invariance of $\L:= \log(\X)$, as follows.  For $a=e^x < b=e^y$,
$|\X \cap (a,b) | = |\L \cap (x,y)|$, with expectation
$\int_x^y \theta \ dt  = \theta(y-x) $ $= \theta \log(b/a) $
$= \int_a^b (\theta/u) \ du$.  Another way to compute this
is to start with the candidate for a scale invariant intensity, $\theta /x \ dx$,
and the
map $x \mapsto c x$, and to apply the change
of variables formula correctly ---  an exercise which
the author usually gets wrong at first!  Note
that scaling ($x \mapsto cx$) the process with intensity
$\theta e^{-x}/x \ dx$ yields the same $\theta$ but a
different exponential decay, i.e. the new intensity is
$\theta e^{-x/c}/x \ dx$, for $0<c<\infty$.

\section{Poisson processes, inversion invariance}
The characterizing property of a Poisson  process is that
for disjoint $I_1,I_2,\ldots,I_k$, the number of points in
$I_j$ for $j=1$ to $k$ are independent, Poisson distributed
random variables.
%PICTURE 1
 This property is preserved by any mapping,
not necessarily one-to-one;  the relevant examples here are
 $x \mapsto x/ \theta$ for $\theta > 0$, $x \mapsto -x$,
$x \mapsto e^x$,  $x \mapsto e^{-x}$,
$x \mapsto \log x$, and  $x \mapsto -\log x$.

Thus in checking that the following four statements about a point
process $\X=\{ X_i \}$ are equivalent, one ingredient
is calculating intensities, and the other ingredient, mapping
the  {\em Poisson} properties, is automatic:
$$
\{X_i \} \mbox{ is scale invariant Poisson } (\theta/x) \ dx
 \mbox{ on } (0,\infty),
$$
$$
\{ \log X_i \} \mbox{ is translation invariant Poisson } (\theta \ dx )
\mbox{ on } (-\infty,\infty),
$$
$$
\{- \log X_i \} \mbox{ is translation invariant Poisson } (\theta \ dx )
\mbox{ on } (-\infty,\infty),
$$
$$
\{1/X_i \} \mbox{ is scale invariant Poisson } (\theta/x) \ dx
 \mbox{ on } (0,\infty).
$$
The above shows how the reflection invariance of a translation
invariant Poisson process on $\BR$ is directly equivalent
to the inversion invariance of a scale invariant Poisson
process on $(0,\infty)$.

\section{Favorite labellings}
For the translation invariant Poisson process, label the points
$L_i$ for $i \in \BZ$ with $L_i<L_{i+1}$ and $L_0 \leq 0 < L_1$.
This seems like the only decent choice.
\begin{equation}\label{picture 2}
- \infty < \cdots < L_{-1}<L_0 \ \leq 0 < \ L_1<L_2<\cdots < \infty \ .
\end{equation}
%PICTURE 2

For the scale invariant Poisson processes, there are two
natural choices for labelling, and it does not seem reasonable to
try to make
one notation fit all situations.

If we focus from one down to zero,
 use $X_i := e^{-L_i}$, so that $X_{i+1}<X_i$,
$X_1$ is the first point to the left of one, and $X_n \downarrow 0$ as
$n \ra \infty$.
\begin{equation}\label{at zero}
 0 < \cdots < X_3 < X_2 < X_1 <1 \ \ \leq X_0 <
X_{-1}< X_{-2}< \cdots < \infty \ .
\end{equation}

%PICTURE 3

If instead we focus from one up to infinity,
use $X_i := e^{L_i}$ so that $X_i < X_{i+1}$ and $X_1$ is
the first point to the right of one:
\begin{equation}\label{at infinity}
 0 < \cdots  < X_{-1} < X_0  \ \ \leq 1 <  X_1 <
X_2< \cdots < \infty \ .
\end{equation}
%PICTURE 4
\section{Spacings  --- for the translation invariant process}
With the labelling of (\ref{picture 2}), the spacings of the points
of the translation invariant process are $L_{i+1}-L_i$ for $i \in \BZ$.
Because of the exceptional behavior of $L_1-L_0$, these spacings
are {\em not} distributed the same as $\ldots, L_{-2}-L_{-1},$
$L_{-1}-L_0,$ $0-L_0$, $L_1-0$, $L_2-L_1,\ldots$ \ , which
are iid, each distributed like an
 exponentially distributed random variable $W$ with
mean $1/ \theta$ and $\p(\theta W>t)=e^{-t}$ for $t>0$.
  Note the familiar ``waiting time for a bus''
paradox: the interval that covers the origin has length $L_1-L_0$
which can be expressed as
the sum of two independent exponentials, $L_1-0$ and $0-L_0$, and
which is equal in distribution to the size biased spacing $W^*$, discussed
in a different context in (\ref{bias}).

Write $W_1:= L_1, W_2:=L_2-L_1, W_3:=L_3-L_2$ etc. so that
for $k=1,2,\ldots$, $ L_k = W_1+\cdots+W_k$.
Recall that ``the (negative)
exponential of an exponential is uniform'', at least
for $\theta=1$.  To check this, note that $t>0 \Leftrightarrow
e^{-t} \in (0,1)$ and $\p(e^{-\theta W_i} < e^{-t}) $
$= \p(\theta W_i > t) $ $=e^{-t}$.  Let
\begin{equation}\label{u from w}
U_i := e^{-W_i} = (e^{-\theta W_i} ) ^{1/\theta} \indist
(\mbox{UNIFORM})^{1/\theta}
\end{equation}
where UNIFORM denotes a random variable uniformly distributed on [0,1].

In lecture, we offered the following closed book QUIZ:  the distribution
of $U_i$  is either Beta$(1,\theta)$ or Beta$(\theta,1)$, but
{\em which one?}

Use the notation (\ref{at zero})
which focusses on the scale invariant process near zero,
i.e.  $X_i =\exp(-L_i)$.  Then for $k=1,2,\ldots$ the expression
\begin{equation}\label{sum of expl}
L_k = W_1+W_2 + \cdots + W_k
\end{equation}
for a sum of independent exponential mean $1/ \theta$ random
variables  is immediately equivalent to
$$
X_k = U_1 U_2 \cdots U_k
$$
with a product of independent (UNIFORM)$^{1/ \theta}$
random variables.

\section{Spacings --- for the scale invariant process}
%PICTURE 5(to show ...$x_2 < x_1 < 1 < x_0 <$ ...)
The process of all spacings for the process
on $(0,\infty)$ with the notation (\ref{at zero}) is defined by
$$
Y_k := X_{k-1} - X_k \in (0,\infty), \ \ \ k \in \BZ.
$$

However, when considering the scale invariant
process restricted to (0,1) or (0,1],
the natural notation gives a different definition for the first
spacing:
\begin{equation}\label{def GEM}
Y_k := X_{k-1}-X_k, \ \ \ k=2,3,\ldots
, \ \ \ \mbox{ but } Y_1 := 1-X_1.
\end{equation}
In terms of the independent $U_1,U_2,\ldots \indist ($UNIFORM)$^{1/ \theta}$,
(\ref{def GEM}) is $Y_1=1-U_1,$
$Y_2 = U_1 -U_1 U_2 = U_1 (1-U_2),$ and in general for $n \geq 1$,
\begin{equation}\label{GEM unif}
Y_n= U_1 U_2 \cdots U_{n-1}(1-U_n).
\end{equation}

\section{Residual allocation, GEM}
Most geneticists have used the notation
$$
Y_1=U_1, Y_2=(1-U_1)U_2, Y_3=(1-U_1)(1-U_2)U_3, \ldots \,
$$
which is the opposite of (\ref{GEM unif}).
We wish to argue that the notation (\ref{GEM unif}),
which respects the sum (\ref{sum of expl}), is preferable.

The product form above, with either notation,  is referred to as a
 ``residual allocation model'', from Halmos 1944 \cite{ram}.
The distribution of $(Y_1,Y_2,\ldots)$ in
(\ref{def GEM}) is called the GEM with parameter $\theta$,
after Griffiths, Engen, and  McCloskey; ---\cite{McCloskey} is the unpublished 1965
thesis by McCloskey; these historical notes are from
chapter
41 by Ewens and Tavar\'e in \cite{encyclopedia}.

The Beta function is defined, for $a,b>0$, by
$B(a,b) := \int_0^1 (1-x)^{a-1}x^{b-1} \ dx$ and
the corresponding distribution, Beta($a,b$), has
density $(1-x)^{a-1}x^{b-1}/B(a,b)$ on (0,1).
Since $ \ B$ is Beta($a,b$) \ if and only if  \ $1-B$ is Beta($b,a$),
we see how  Beta(1,$\theta$) and
Beta($\theta$,1) can easily be confused.

\section{The Poisson-Dirichlet process}
We write $(V_1,V_2,\ldots)$ for the Poisson-Dirichlet process
with parameter $\theta$.
Survey references for the Poisson-Dirichlet process include
\cite{pityor97,abt}, and Chapter 41 in \cite{encyclopedia}.
We do not assume the reader knows
this process; we will give several characterizations of it, including
one that relates it to the scale invariant Poisson process with parameter
$\theta$, following a review of its history.
When the Poisson-Dirichlet process was first studied, none of its
connections to the scale invariant Poisson process was
explicitly noted.

One might view the Poisson-Dirichlet process with $\theta=1$
as implicit in the 1930 paper  of Dickman \cite{dickman}
describing the limit distribution of the largest prime
factor of a randomly chosen integer.  However if one insists
on an explicit description of the process, it seems that the
first appearance of the Poisson-Dirichlet
is in Billingsley 1972 \cite{billing large}, which proves that
\begin{equation}\label{billingsley}
\ \ \left( \frac{\log P_1}{\log n } , \frac{\log P_2}{\log n}, \ldots \right) \todist
(V_1,V_2,\ldots),
\end{equation}
where the limit is the Poisson-Dirichlet process with $\theta=1$.
Here, $P_i$ is the $i^{th}$ largest prime factor
of an integer $N$ chosen uniformly from 1 to $n$, using $P_i=1$ when
$i$ is greater than the number $\Omega(N)$ of prime factors including
multiplicities. (For example, if the random integer is 12, then $P_1=3,$
$P_2=P_3=2, 1=P_4=P_5=\cdots \ . $)
Billingsley specified the limit process in terms of
the joint density of $(1/V_1,1/V_2,\ldots,1/V_k)$.
The second published proof of Billingsley's result, \cite{DG},
is based on a size-biased permutation, and is especially robust, leading
to a proof \cite{primecon}
of  the analogous Poisson-Dirichlet convergence
for general $\theta \neq 1$, in which the random integer is conditioned
 on the large deviation that the
number of distinct prime factors is $\theta \log \log n$.

Ferguson \cite{ferguson} in
1973 described a class of processes that includes what we now know as
the Poisson-Dirichlet; a recent survey of this is Pitman\cite{pitman-blackwell}.

Kingman \cite{kingman75} in 1975 and Ignatov \cite{ignatov} in
1982  gave the first direct connection
between the Poisson-Dirichlet process and the scale invariant
Poisson process: they can be coupled, with $V_i =$ the
$i^{th}$ largest of the spacings $Y_1,Y_2,\ldots$ as in (\ref{def GEM}),
i.e.
\begin{equation}\label{rank gem}
  (V_1,V_2, \ldots) \indist \mbox{ RANK}(1-X_1,X_1-X_2,X_2-X_3,\ldots).
\end{equation}
See \cite{pityor96} for generalizations.

Vershik and Schmidt \cite{vershik1,vershik2} in 1977  showed that
$$
(\theta=1) \ \ \ \ \  \left( \frac{L_1}{ n } , \frac{L_2}{n}, \ldots \right) \todist
(V_1,V_2,\ldots)
$$
where $L_i$ is the length of the $i^{th}$ longest cycle of
a random permutation of $n$ objects, choosing with all
$n!$ possibilities equilikely.

Aldous in 1983 \cite{aldous mapping}
 gave the analogous Poisson-Dirichlet limit for
the component sizes in a random  mapping on $n$ points, with
all $n^n$ maps equilikely; here $\theta=\frac{1}{2}$.

Hansen \cite{hansen} gave a general treatment of decomposable
combinatorial structures having a Poisson-Dirichlet limit; see
also \cite{abtlocal}.

The following comes from the 1993 book \cite{kingman},
``Poisson processes,''
by Kingman.   In 1968 Moran considered the Poisson process
with intensity $\theta \exp(-x)/x \ dx$ on $(0,\infty)$.  The points of this
process can be labelled $\sigma_i$ for $i=1,2,\ldots$ and
$0< \cdots < \sigma_2 < \sigma_1$.  The sum $\sigma :=
\sigma_1+\sigma_2+\cdots$ has a Gamma distribution with
parameter $\theta$, and is independent of the rescaled
vector $(\sigma_1/ \sigma, \sigma_2 / \sigma,\ldots)$,
which has the Poisson-Dirichlet distribution.   There is a somewhat similar
model in physics, discussed in \cite{ruelle}, 227-228, with a Poisson
process on $(0,\infty)$ having intensity $c x^{-c-1} \ dx$ for
a constant $c \in (0,1)$, instead of $\theta \exp(-x)/x \ dx$.  It
is still the case that a.s. the points can be labelled
$\sigma_i$ for $i=1,2,\ldots$ with
$0< \cdots < \sigma_2 < \sigma_1$, and  $\sigma :=
\sigma_1+\sigma_2+\cdots< \infty$,
and one may form the rescaled
vector $(\sigma_1/ \sigma, \sigma_2 / \sigma,\ldots)$.
However, in this model, $\e \sigma = \infty$, in contrast to
$\e \sigma = \theta$ for the Moran model.
A unified treatment which includes the two models
is \cite{pitman multiple}.

\section{Joint density for the Poisson-Dirichlet}
The joint density $f_k$ of the first $k$ coordinates
$(V_1,V_2,\ldots,V_k)$ of the Poisson-Dirichlet distribution
has support $\{ (x_1,\ldots,x_k): \
1>x_1> \cdots > x_k > 0$ and $x_1+\cdots+x_k <1 \}$.
At such points, for the case $\theta=1$,
\begin{equation}\label{pd 1}
f_k(x_1,\cdots,x_k) = \rho \left( \frac{1-x_1-\cdots - x_k}{ x_k }
 \right)  \frac{ 1 }{ x_1x_2\cdots x_k }
\end{equation}
where $\rho$ is Dickman's  function
\cite{dickman,tenenbaum},
%PICTURE 6
characterized by
$\rho=0$ on $(-\infty,0)$, $\rho=1$ on [0,1],
$\rho$ continuous on $(0,\infty)$, and
$\rho'(u)=-\rho(u-1)/u$ for $u>1$.

The area under the graph of $\rho$ is $e^\gamma$,
where $\gamma$ is Euler's constant,
and $\rho \geq 0$ everywhere, so $g(u):=
e^{-\gamma} \rho(u)$ is a probability density, which
might naturally be called the Dickman distribution.

\section{What is the ``Dickman distribution''?}
Dickman in 1930 showed that $\log P_1 / \ \log n \todist V_1$,
where $\p(V_1 < 1/u)=\rho(u)$, i.e.~the Dickman function
$\rho$ gives the tail probabilities for $1/V_1$.  So if you hear someone
refer to a random variable having ``the Dickman distribution,''
you may be confident that the speaker refers to
 one of the following three, but which one?

\begin{itemize}

\item the random variable $T$ with density $g(u) := e^{-\gamma} \rho(u)$
on $(0,\infty)$

\item $V_1$, with density $\rho((1-x)/x)/x$ on $(0,1)$

\item $1/V_1$, with density $\rho(u-1)/u$ on $(1,\infty)$

\end{itemize}
The Dickman function decays superexponentially fast.
Hildebrand 1990 \cite{hildebrand}
gives the following asymptotic expansion
of the Dickman function.  Write $L := \log u$ and
$M:= \log \log u$. As $u \ra \infty$,
$$
\frac{-\log \rho(u)}{ u} = L  + M - 1 + \frac{M}{  L}
- \frac{1}{ L} - \frac{M^2}{2L^2} + \frac{M}{L^2}
- \frac{2}{L^2} + O \left( \frac{M^2}{L^3} \right),
$$
so that $\rho(u) = \exp(-u \log u - u \log \log u + u - \cdots) $.

\section{A key random variable: \ T, the sum of the points in (0,1)}

%PICTURE 7, OF MONKEY

Let $T$ be the sum of the locations of all points of the
scale invariant Poisson process, with intensity $\theta / x \ dx$,
restricted to $(0,1)$.  For the case $\theta=1$, this is the random variable
in the previous section, with density $e^{-\gamma} \rho(u)$.  For
any $\theta>0$,
$$
T:= X_1 + X_2 + X_3 + \cdots
$$
$$
= U_1 + U_1 U_2 + U_1 U_2 U_3 + \cdots
$$
$$
=U_1 \ (1 + U_2 + U_2 U_3 + \cdots)
$$
$$
\indist U \ (1+ T')
$$
with $T', U$ independent, $T' \indist T$, and
$U \indist ($ UNIFORM $)^{1/ \theta}$.
Thus it is elementary to get an integral equation
involving the density $g$ of $T$, although
this equation is not especially tractable.

Another approach is via size biasing; the following is
taken from \cite{abt}.  In general,
if $X$ is a non-negative random variable with mean
$\mu \in (0,\infty)$, then the size biased random
variable $X^*$  is characterized by
\begin{equation}\label{bias}
\e h(X^*)=\e(X h(X))/\mu;
\end{equation}
if further $X$ has density $g$ then $X^*$ has density
$x g(x) / \mu$.  If
$Z$ is a Poisson distributed random variable then $Z^* \indist 1+Z$ and
for any $x>0$, $(xZ)^* \indist x + (xZ)$.
If $Z$ is a sum of independent nonnegative random variables, with
$\e Z < \infty$, then $Z^*$ is likewise a sum of independent random variables,
using the same summands except that one summand is size biased, and the
choice of which summand to bias is made with probability proportional to its
contribution to $\mu$.  Thus if $T$ is the sum of locations of points in a
Poisson process on $(0,\infty)$ with intensity $f(x) \ dx$,
such that $\mu := \e T =
\int x f(x) \ dx \ < \infty$, then $T^*$
is formed by choosing a location $x$ with probability $x f(x) \ dx / \mu$
where the summand is size biased by deterministically adding in $x$. Hence
the density of $T^*$ at $t$ is
$$
\frac{t g(t)}{\mu} = \int_0^\infty \frac{x f(x) \ dx}{\mu} \ g(t-x)
$$
For the case of the scale invariant Poisson process restricted to
(0,1) we have $\e T = \int_0^1 x (\theta/x) \ dx \ = \theta$, and
$$
\frac{t g(t)}{\theta} = \int_0^1 \frac{x (\theta/x) \ dx}{\theta} \ g(t-x) \ = \int^t_{t-1} g(u) du,\ \  \ \ t>0,
$$
Hence the differential-difference equation satisfied by $g$ is
$$
t g^\prime(t) + (1-\theta) g(t) + \theta g(t-1) =
0,\ \  \ \ t>0.
$$

Another approach is that since $T$ is the sum of locations
of the Poisson process with intensity function
$f(x) = \bone(0<x<1) \theta / x$,
it has Laplace transform
$$
\e \exp(-sT) = \exp\left(- \int_0^\infty (1-e^{-sx}) f(x) \ dx  \right)=
\exp\left( - \theta \int_0^1 \frac{1 - e^{-sx} }{x } \ dx \right)
$$
This leads to properties of the density $g$; see Vervaat 1972
\cite{vervaat},  Watterson \cite{watterson}, and Hensley \cite{hensley}.

\section{The scale invariant process of sums, and class L}
For any fixed $\theta>0$,
our favorite random variable $T$, which is the sum of locations of
points in (0,1).  It satisfies $T \indist T_1 \indist (1/t) T_t$,  where
for $t > 0$, $T_t$ is defined as the sum of locations of
all points in $(0,t]$, for the scale invariant Poisson process with
intensity $\theta/x \ dx$. Of course we take $T_0 \equiv 0$.
As a process $(T_t)_{t \geq 0}$
\begin{itemize}

\item  increases by jumps

\item  at time $t$, can only stay constant or jump up by $t$

\item has a jump occurring in $(t,t+dt)$ with probability
$(\theta/t) \ dt \ ( 1+o(1))$

\item  has independent increments

\item is self-similar with index 1: for $t>0$,
  $(T_{ts})_{s \geq 0} \indist t^1 (T_s)_{s \geq 0}$

\end{itemize}
In particular, writing $g_t$ for the density of $T_t$ and $g$ for
the density of $T$, from $T_t \indist t \, T$ \ we have, for any
$t>0$,
$$
g_t(x) = \frac{1}{t} \ g(\frac{x}{t}).
$$

From Feller II \cite{feller2}, the L\'evy class L of infinitely
divisible distributions consists of those which are limit
distributions of $\{ S_n^*\}$ where $S_n := X_1+\cdots +X_n$,
with independent (not necessarily identically distributed) $X_1,X_2,\ldots$,
and $S_n^*=(S_n~-~b_n)~/~a_n$ for constants $a_1,a_2,\ldots$ and $b_1,b_2,\ldots$
such that $a_n \ra \infty$, $a_{n+1}/a_n~\ra~1$.
To see that for each $\theta>0$ the random variable
$T$ is in class L, take $X_n := T_n-T_{n-1}$, so that $S_n=T_n$, and
take $b_n =0, a_n=n$ so that $S_n^*:=T_n/n \indist T$ for all $n$.
We thank Larry Shepp for an enjoyable conversation on this topic.

We cannot understand a technical/aesthetic issue involving the
process $(T_t)$:  why does the natural
random Stieltjes measure $dT_{(\cdot)}$ always lose out
to the competing counting measure?  That is, to
encode the scale invariant Poisson process as a random
measure,  the usual choice is  the
counting measure $\sum_{i \in \BZ}
\delta_{X_i}(\cdot)$.  This has to be taken as a random measure on
$(0,\infty)$ with the point at zero not in the underlying space; there is
infinite mass in every neighborhood of zero.
Another sensible choice would  be
$d{T_{(\cdot)}} \equiv \sum_{i \in \BZ} X_i
\ \delta_{X_i}(\cdot)$, which gives a sigma-finite
random measure on $[0,\infty)$, with finite mass near zero.

\section{Conditioning on $T=s$ in general}

Consider {\em any} Poisson process on (0,1], having intensity
$f(x) \ dx$, such that $T_t $, defined to be the sum of the locations of all
points in $(0,t)$, has a density $g_t$.  Write $T := T_1$  and $g :=g_1$.  Assume that $g$ is strictly positive on
$(0,\infty)$.  Label the points in (0,1) so that
$1>X_1>X_2>\cdots > 0$.

The joint density of $X_1,\ldots,X_k$ is
$$
f(x_1) f(x_2) \cdots f(x_k) \exp \left( - \int_{x_k}^1 f(u) \ du \right),
$$
supported by points
in $(0,1)^k$ where $x_1>x_2\cdots > x_k$.  The first $k$
factors correspond to requiring points at $x_1,\ldots,x_k$, and the
last factor corresponds to demanding no other points in
$(x_k,1)$.
Thus for any $s>0$
 the joint density of $X_1,\ldots,X_k$ conditional on $T=s$ is
\begin{equation}\label{general cond}
\frac{g_{x_k}(s-x_1-\cdots-x_k)}{g(s)} \
f(x_1) f(x_2) \cdots f(x_k) \exp \left( - \int_{x_k}^1 f(u) \ du \right),
\end{equation}
supported at $1>x_1>\cdots>x_k>0$ with $x_1+\cdots+x_k <s$.
For the special case of the scale invariant Poisson process with
$f(x) = \theta / x$, this is
$$
\frac{1}{x_k} \ \frac{g \left( (s-x_1-\cdots - x_k)/x_k \right) }
{g(s)}
 \ \frac{ \theta }{x_1} \ \frac{\theta}{x_2} \ \cdots
\frac{\theta}{x_k}
e^{- \theta \log(1/x_k)}
$$

\section{Conditioning on $T=1$ for the scale invariant process}
For the special case $s=1$ the conditional joint density of
$X_1,\ldots,X_k$ given $T=1$ simplifies to
$$
\frac{g \left( (1-x_1-\cdots - x_k)/x_k \right) }
{g(1)}
 \
\frac{\theta^k}{x_1 x_2 \cdots x_k}
x_k^{\theta -1}
$$
\begin{equation}\label{pd general}
= g\left(\frac{1-x_1- \cdots - x_k}{x_k}\right) \,
\frac{e^{\gamma \theta} \,
\theta^k \, \Gamma(\theta) \, x_k^{\theta-1}}  {x_1 x_2 \cdots x_k},
\end{equation}
which is the joint density of the Poisson-Dirichlet; the final equality
is just giving an explicit formula for the
normalizing constant $g(1)$, which for $\theta \neq 1$ comes from
\cite{watterson}.  In the special case $\theta=1$, (\ref{pd general})
simplifies to (\ref{pd 1}).

In summary,
if one starts from the joint density of the Poisson-Dirichlet, with a
differential-difference equation or Laplace transform to characterize
the difficult factor, then it is very easy to prove that the
Poisson-Dirichlet process is the scale-invariant Poisson process, restricted to
(0,1), and conditioned on the event $T=1$:  for any $\theta > 0$
\begin{equation}\label{t1}
(V_1,V_2,\ldots) \indist ( \ (X_1,X_2,\ldots) \ | \ T=1 \ ), \ \mbox{ i.e.} \ \rm{PD}(\theta) \indist (\rm{PP}(\theta) \ |  \ T=1 \ ).
\end{equation}

The above result first was written in the 1996 version of \cite{abt}.
One motivation for considering it is the analogy to the result that
\begin{equation}\label{cond permute}
(C_1(n),\ldots,C_n(n)) \indist ( \ (Z_1,\ldots,Z_n) \ | T_n = n \ ),
\end{equation}
relating the joint distribution of cycle counts for a uniformly
distributed random permutation, to the joint distribution of
$Z_1,Z_2,\ldots$, a sequence of independent Poisson random variables
with $\e Z_i = 1/i$, conditioned on $T_n=n$, where $T_n := Z_1 + 2 Z_2 +
\cdots + n Z_n$.

Just as (\ref{t1}) is a continuum analog of
(\ref{cond permute}),  the Moran process representation implies
\begin{equation}\label{moran}
(V_1,V_2,\ldots) \indist ( \ (\sigma_1,\sigma_2,\ldots) \ | \ \sigma=1 \
),
\end{equation}
which
is a continuum analog of the relation exploited by
Shepp and Lloyd \cite{shepp and lloyd} in 1966 to
analyze random permutations,
\begin{equation}\label{shepp and lloyd}
(C_1(n),\ldots,C_n(n),0,0,\ldots)
 \indist ( \ (Z_1,Z_2,\ldots)  \ |\  T_\infty = n \ ).
\end{equation}
Here $T_\infty=Z_1+2Z_2+\cdots$,
and the $Z_i$ are independent Poisson with $\e Z_i = z^i/i$;
this is valid for any $z \in (0,1)$ but especially useful with $z=1-1/n$.
Furthermore, (\ref{t1}) is to (\ref{moran}) as
(\ref{cond permute}) is to (\ref{shepp and lloyd}).
Loosely speaking,  (\ref{moran}) and (\ref{shepp and lloyd})
both avoid the divergence of $\int^\infty x^{-1}$ by throwing
in an exponentially decaying factor which does not affect the conditional
distribution, while (\ref{t1}) and (\ref{cond permute}) avoid the
divergence of $\int^\infty x^{-1}$ by not going out to infinity.

Another elementary proof of  (\ref{t1}) can be given by comparison to
the Moran representation; this is done in \cite{T1}.  The general
conditioning formula (\ref{general cond}) is also used in \cite{perman}.

\section{More discrete analogs}
The role of the scale invariant Poisson process in number theory may be
found explicitly in
DeKoninck and Galambos \cite{k+galambos},
with a theorem showing that the process of logarithms of  the intermediate
prime divisors of a random integer chosen uniformly from 1 to $n$
converges in distribution to the scale invariant Poisson process, with
$\theta=1$.  An attempt to metrize this result \cite{primedw}
leads to the following,
which shows how  the scale invariant
Poisson process with $\theta=1$ is very close to
the discrete limit process that arises in number theory,
with the number $Z_p$ of occurrences of each prime $p$ being independent, with
$\p(Z_p=k) = (1-1/p) p^{-k}$.  Theorem: it is possible
to couple random values $Q_i$ for $i \in \BZ$, and the scale invariant
Poisson process, so that
$$
\e \sum_{i \in \BZ} | X_i - \log Q_i | \ < \infty
$$
with each $Q_i$ either one or prime, and for each $p$, $Z_p=\sum_i
\bone(Q_i=p)$.

It is not obvious how the usual model of number theory,
picking an integer uniformly from 1 to $n$, is related to the
independent, geometrically distributed $Z_p$ by
something like conditioning; this is explained in \cite{AMS1}.
The overall theme there is that discrete ``logarithmic'' combinatorial
structures, such as random permutations, random mappings, and random
polynomials over finite fields, together with prime factorizations
of uniformly chosen random integers,
have a dependent component size counting
process, that is derived from an independent
process by conditioning or something analogous.
The discrete dependent processes,
rescaled, are close to the Poisson-Dirichlet processes.  The discrete
independent processes, rescaled, are close to the scale invariant
Poisson processes.  For both the discrete and continuous situations, the
dependent processes are obtained from the independent processes by
conditioning, or something like conditioning --- this is the point of
studying (\ref{t1}).

\section{Dependent versus independent: total variation distance}
Recall, the dependent process is the Poisson-Dirichlet
$(V_1,V_2,\ldots)$ with $V_1+V_2+\cdots=1$, and the related
independent process
is the scale invariant Poisson, restricted to (0,1), with $X_1+X_2+\cdots = T$,
a random variable having a strictly positive density on $(0,\infty)$.
For $\beta \in [0,1]$ let
$$
H_\theta(\beta) := d_{TV}(\ \{V_1,V_2,\ldots \} \cap [0,\beta], \
  \{ X_1,X_2,\ldots \} \cap [0,\beta] \ ).
$$
Note that $H_\theta(0) = d_{TV}( \emptyset, \emptyset) = 0$, and
$$
H_\theta(1) = d_{TV}(\ (V_1,V_2,\ldots), (X_1,X_2,\ldots)\  ) = 1
$$
because $1=\p(V_1+V_2+\cdots=1)$ and
$0=\p(X_1+X_2+\cdots=1)$.
For the case $\theta=1$
$$ H_1(1/ \, 1.9) = .4968\ldots \ ,  H_1(1/ \, 2) = .4454\ldots \ ,
 H_1(1/ \, 3) = .1114\ldots \ ,
$$
$$
 H_1(1/ \, 3.5) = .0471\ldots \ ,  H_1(1/ \, 4) = .0184\ldots \ .
$$
Informally, $H_1(.25)=.0184\ldots$ means that
if you are shown one sample of
a process which is either the Poisson-Dirichlet
($\theta=1$) or else the scale invariant Poisson ($\theta=1$),  but you
only get to observe the set of components of size at most .25, then
your edge over the house is at most 1.84 percent:  if the examiner
picks his distribution using a fair coin, then you have at most
a 50.92 percent chance of answering correctly.

There is a phase transition in the qualitative
behavior of $H_\theta(0+)$, with linear decay
for $\theta \neq 1$, and superexponential decay for
$\theta =1$.  In detail, from \cite{ast}, p.~1368, combined
with \cite{T1,stark,AS}, as $\beta \ra 0+$,
\begin{eqnarray}\label{linear}
  \frac{1}{\beta} H_\theta(\beta) & \ra & |1-\theta|\ \frac{e^{-\gamma \theta}}{\Gamma(\theta)} \ \frac{\theta^\theta}{1+\theta} >0 \mbox{ for } \theta \neq 1,
\\
\label{super}
 - \beta \log H_1(\theta) & \sim & \log(1/\beta) \ra \infty.
\end{eqnarray}

\section{Invariance principle for total variation distance}
The function $H_\theta(\beta)$ gives the total variation
distance between the dependent and independent processes,
i.e.~ the Poisson-Dirichlet process with parameter $\theta$, and
the scale invariant Poisson process with parameter $\theta$,
when observing components of size at most $\beta$.
These functions $H_\theta(\cdot)$, which can be most
easily defined in terms of the scale invariant Poisson,
showed up first in combinatorics.
Consider random permutations of $n$ objects, observing
cycles of length at most $\beta n$: write $C_i(n)$ for
the number of cycles of length $i$, and $Z_i$ for a
random variable which is Poisson, mean $1/i$, with
$Z_1,Z_2,\ldots$ independent.  Then \cite{tv1,stark} for any fixed
$\beta \in [0,1]$, as $n \ra \infty$,
$$
d_{TV} (\  (C_1(n),C_2(n),\ldots, C_{\lfloor \beta n \rfloor}(n)), \
(Z_1,Z_2,\ldots,Z_{\lfloor \beta n \rfloor}) \ ) \ \ra H_1(\beta).
$$
Similarly \cite{stark} with $\theta=\frac{1}{2}$,
for the component counts for a random mapping, comparing
to the independent limit process, and observing components of size at most
$\beta n$,
$$
d_{TV} \ra H_{\frac{1}{2}}(\beta),
$$
and also similarly with $\theta=1$ for the factorization of a random
polynomial of degree $n$ over a finite field, observing
factors of degree at most $\beta n$.

For the factorization of an integer chosen uniformly from 1 to $n$ the
analogous result holds \cite{AS,sieve}, again with $\theta=1$.
Here, we observe prime factors, with multiplicity, jointly for
all primes $p \leq n^\beta$.  The independent comparison process has
$Z_p$ which are geometrically distributed with $\p (Z_p \geq k) =
p^{-k}$.

\section{Buchstab and the explicit formula for $d_{TV}$}

The possiblity of numerically evaluating the limiting
total variation distance $H_\theta(\beta)$
arises thanks to relation (\ref{t1}),
which makes it possible to equate $H_\theta(\beta)$,
defined as the total variation distance between
processes, with the total variation distance
between two random variables.  Recall that the
process $(T_s)_{s \geq 0}$ has independent
increments, so that the $T \equiv T_1 $ in conditioning
on $T=1$ can be expressed as the sum
$T=T_\beta + (T-T_\beta)$ with two independent
summands.  This leads to
$$
H_\theta(\beta) := d_{TV} \mbox{ for processes }
$$
$$ = d_{TV}
\mbox{ for random variables } \equiv d_{TV}( T_\beta, ( T_\beta \ | \ T=1 \ ) \ )
$$
Manipulation of the densities of the independent
summands leads
to the explicit expression for $H_\theta(\beta)$, which for $\theta=1$
can be expressed as:
\begin{eqnarray*}
2H_1(\beta) & = & e^\gamma \e | \omega(u-T)-e^{-\gamma}| \ + \rho(u) \\
 & = & \int_{t>0} |\omega(u-t)-e^{-\gamma}| \rho(t) \ dt \ \ + \rho(u).
\end{eqnarray*}
Here, $\rho$ is Dickman's function, and
and $\omega$ is Buchstab's function, with
$\omega$ continuous on $(1,\infty)$, with
$\left(u \omega(u)\right)'=\omega(u-1)$ for $u>2$ and $u\omega(u)=1$ for
$1 \leq u \leq 2$; see e.g.
\cite{tenenbaum}.  Both functions
were described some sixty years ago in number theory,
and each can also be described in terms of the
scale invariant Poisson process for $\theta=1$:
$\rho$ is $e^\gamma$ times the density function
of $T$, and for $u > 1$ with
 $\beta := 1/u$,  $\omega(u)$ equals the density of
$T-T_\beta$, evaluated at $1-$, i.e.
$\omega(u)=\lim_{\Delta t \ra 0+} \p( T-T_{1/u} \in (1-\Delta t,1) \, ) / \Delta t$.
To state the connection with number theory,
$\Psi(x,y)$ counts positive integers less than or equal to $x$ whose largest
prime factor is less than or equal to $y$, $\Phi(x,y)$ counts positive integers
less than or equal to $x$ with no prime factor less than or equal to $y$,
and for $u > 1$ as $n \ra \infty$,
\begin{equation}\label{dick buch}
\frac{1}{n} \ \Psi \left(n,n^{1/u}\right) \ra \rho(u), \ \
\frac{\log n}{u} \ \frac{1}{n} \  \Phi \left(n,n^{1/u} \right) \ra \omega(u).
\end{equation}

For every $\theta$, $H_\theta$ is continuous and strictly monotone, mapping
[0,1] into itself, with $H_\theta(0)=0$ and $H_\theta(1)=1$, but only for $\theta=1$ is
it the case that $H_\theta(\beta)=o(\beta)$ as $\beta \downarrow 0$.

From $H_\theta(\beta)\ra 0$ as $\beta \downarrow 0$ it follows that
the Poisson-Dirichlet process, ``blown up'', coverges in distribution
to the scale invariant Poisson.  That is, with $\V := \{ V_1,V_2,\ldots \} \subset (0,1)$,
and $\X:= \{ X_i: \ i \in \BZ \} \subset (0,\infty), $
as $v \ra \infty$, $v \V \todist \X$.  For a proof: for any fixed $x > 0$, for
$v \geq x$ we have
\begin{equation}\label{blow up}
d_{TV}(\ v \V \cap (0,x), \X \cap (0,x) \ ) = d_{TV}(\  \V \cap (0,x/v), (v^{-1}\X) \cap (0,x/v) \ )
\end{equation}
$$
=H_\theta(x/v).
$$
where we use the scale invariance of $\X$ to see that last equality. Thus
for fixed $x$, $d_{TV}(v \V \cap (0,x), \X \cap (0,x) \ )  \ra 0$
as  $v \ra \infty$, which is quite a bit
stronger than distributional convergence.

\section{Insertion---deletion distance}
That $H_\theta(1)=1$ says that by looking at a single sample
one can a.s.~ tell the random set
$\V := \{ V_1,V_2,\ldots \}$ from $\Y := \X \cap (0,1) =
\{ X_1, X_2,\ldots \}$.
Could the two random sets nevertheless be close?
Consider  \cite{primedw} the
Wasserstein metric $d_W$ using counts of insertions and
deletions needed to convert one set to the other:
$$
d_W := \min \ \e | \V \triangle \Y  |
$$
$$ = \min ( \ \e|\V\setminus \Y| + \e | \Y \setminus \V| \ ).
$$
The minimum is taken over all couplings of the Poisson-Dirichlet process and the
scale invariant process, both with parameter $\theta$.

For $\theta=1$, this can be fully understood:
$d_W = 2$, and there is a coupling such that
$$
 \V \setminus \{ V_J \}
= \Y   \setminus
 \{ X_{I_1},\ldots, X_{I_D} \ \}
$$
where the number $D$ of points of the scale invariant Poisson
process restricted to (0,1) that need to be deleted is
random, with $\e D =1$.

That $d_W \geq 2$ for
$\theta=1$ starts with a proof that a.s.~at least one deletion from
$\V$ is needed,  because with probability one no subsum of the points of
$\X$ has the value $V_1+V_2+\cdots=1$.  This argument requires only that $\theta \leq
1/ \, \log 2$; see problem \ref{problem log 2}.  For the other
half of the lower bound, $\theta=1$ is needed, so
that $\X$ and $\V$ have the same intensity, and hence the one
required deletion from $\V$ must be matched by deletions from $\X$, averaging
one in number.

The coupling
showing that  $d_W \leq 2$ for $\theta=1$
can be given explicitly. In this coupling,
the deleted compenent $V_J$ is simply
the first pick from $\V$ in a size biased permutation, i.e. given the
value of $(V_1,V_2,\dots)$, the conditional probability that
$V_J=V_k$ is $V_k$.   The entire coupling
 is the continuum analog of the
Feller coupling for random permutations \cite{tv3},
 and depends on the ``scale
invariant spacing lemma''.

\section{Scale invariant spacing lemma}
Label the scale invariant Poisson process as in (\ref{at infinity}),
with $X_i < X_{i+1}$ for $i
\in \BZ$, and let $Y_i := X_{i+1} - X_i$.
$$
\mbox{ LEMMA.  For any } \theta > 0, \ \ \{ X_i: \ i \in \BZ \ \}
\indist \{ Y_i: \ i \in \BZ \ \}
$$
and a.s. all the $Y_i$ are distinct.

This is the continuous analog  of a property of the Feller coupling:
  if $\xi_1,\xi_2,\ldots$ are independent Bernoulli with
$\p(\xi_i=1)=\theta/(\theta+i-1)$, and $Z_k$ is defined as the number of
$k$-spacings between consecutive ones in $\xi_1,\xi_2,\ldots$, then
the $Z_1,Z_2,\ldots$ are independent, Poisson, with $\e Z_k = \theta/k$.
In both the discrete and continuous cases, one starts with an
independent process, and the surprise is that the spacings also form an
independent process.
In the discrete case,  while each of the
processes $\xi_1,\xi_2,\ldots$ and $Z_1,Z_2,\ldots$
has independent coordinates, the two processes are not the same.

The proof \cite{primedw}
of the scale invariant spacing lemma is based on
the Poisson process on $(0,\infty)^2$ with intensity
$\theta e^{-wy} \ dw \ dy$, whose projections on
each coordinate axis give realizations of the scale
invariant Poisson process.  Labelling the points in
decreasing order of their $w$  coordinates gives
the $y$ coordinates in a permutation such that
$X_n := \sum_{i \leq n} Y_i$ forms a simple
point process $\X = \{ X_n\!: n \in \BZ \}$
on $(0,\infty)$ whose spacings,
by construction,  are the
points $Y_i$ of the scale invariant Poisson process.  It is
then  a calculation that the distribution of  $\X$ is
also scale invariant Poisson.  Since the joint density
function of the $(W_i,Y_i)$
can be written $\theta e^{-wy}$
$= (\theta / y) \ ( y e^{-wy}) $ ($=f_Y(y)\  f_{W|Y}(w)$,)
 which says
that the conditional distribution of the $w$ coordinate, given
$y$, is exponential with mean $1/y$, our construction
may be viewed as using a size-biased permutation of $\{ Y_i \}$;
see \cite{ppy}.

\section{Problem session: what other processes are their own spacings?}

%\noindent {\bf Open problem 1}
\begin{problem}
What other intensity functions
 $f(x)$ on $(0,\infty)$ beside
those of the form $f(x)=\theta/ x$ lead to Poisson processes whose
spacings are again Poisson processes?
\end{problem}

%\vskip .05in
Tom Kurtz, upon being asked the above question, which is still open, immediately
countered with a very different question:

\begin{problem}\label{prob 2}
What point processes on $(0,\infty)$ are equal in distribution to their
own spacings?  Formally, what point processes with points $X_i$ for $i
\in \BZ$ with $X_i < X_{i+1}$ for all $i$ have
\begin{equation}\label{random spacing}
  Y_i := X_{i+1} - X_i \mbox{ are all distinct, a.s., and }
 \{ X_i: \ i \in \BZ \ \}
\indist \{ Y_i: \ i \in \BZ \ \}.
\end{equation}
\end{problem}

The above question is reminiscent of an intriguing, still open
1978 question from Liggett \cite{Liggett}, about invariant random measures for independent particle systems,
a conjecture that $MP=M$ in distribution implies the existence of
a nontrivial solution of $MP=M$ almost surely.

The answer to problem \ref{prob 2} involves mixtures of {\em distributional}
solutions, such as the scale invariant Poisson processes and perhaps some
others, and {\em deterministic}
solutions.
During a problem session at the DIMACS workshop, I
asked what are all simple deterministic solutions, i.e. doubly infinite
sequences, such as $x_i:= 2^i$,
 with $0 < \cdots <x_{-1}<x_0<x_1<x_2<\cdots < \infty$
such that
\begin{equation}\label{simple det}
y_i := x_{i+1}-x_i  \mbox{ are all distinct,  and }\{y_i: \ i \in \BZ \} =
\{x_i: \ i \in \BZ \}.
\end{equation}
G\'abor Tardos and L\'aszl\'o Lov\'asz from the
audience quickly contributed a partial solution: for any fixed
$k \geq 0$ take
\begin{equation}\label{geometric}
x_i := b^i,  \ \ \mbox{ where }  b>1  \mbox{ solves }b^{k+1}-b^{k}=1,
\end{equation}
so that
$b$ is 2 when  $k=0$, and $b $ is the golden ratio when $k=1$.
They observed that $y_i := x_{i+1}-x_i<x_{i+1}$ implies $y_i = x_{i+1}-
x_i \leq x_i$, i.e. $x_{i+1} \leq 2 x_i$, so that the example with
$x_i=2^i$ is extreme.
They also showed how ``entrance solutions'' may be extended: working from left to
right, the next spacing is chosen from the previously defined $x_k$,
excluding those already used as spacings, i.e.~$y_i \in \{
\ldots,x_{i-2},x_{i-1},x_i \}
\setminus \{ \ldots, y_{i-2},y_{i-1} \}$.  For example, $b=(1+\sqrt{5})/2$
and $x_i=b^i$ for $i \leq 0$ is an entrance solution with $y_i = x_{i-1}$
for $i<0$, which may be extended with an infinite series of two-way choices:
first $x_1=x_0+y_0$ where $y_0 \in \{ x_0,x_{-1} \}$, then $x_2=x_1+y_1$ where
$y_1 \in \{ x_1,x_0,x_{-1} \} \setminus \{ y_0 \}$, and so on.
To describe their idea in more detail, focus on the
implied permutation, as follows.

Any  solution  of (\ref{simple det}) determines a permutation $\pi$ of $\BZ$,
such that $ \forall i \in \BZ$, $y_i = x_{\pi(i)}$ , i.e.  $\forall i \in \BZ$
\begin{equation}\label{pi}
x_{i+1} = x_i + x_{\pi(i)}  \in (0,\infty).
\end{equation}
   Note that always $\pi(i) \leq i$ since
$y_i \leq x_i$. The geometric solutions in (\ref{geometric})
correspond to permutations $\pi$ with
$\pi(i)=i-k$.
In constructing deterministic solutions, for any given $\pi$ with
\begin{equation}\label{perm}
\pi \mbox{ permutes } \BZ, \ \mbox{ and } \ \pi(i) \leq i \ \ \ \forall i \in \BZ,
\end{equation}
{\em if} there is an ``entrance solution'' from zero, then the solution can be uniquely extended out to infinity.  That is, if for some $k$
there are $0 < \cdots < x_{k-2}<x_{k-1}<x_k$ satisfying  (\ref{pi}) for all
$i<k$, then recursively defining  $x_{k+1}, x_{k+2},\ldots$ by (\ref{pi}) for $i=k,k+1,\ldots$ involves no
problem --- the solution cannot explode to plus infinity in finite time.  Running this
argument in the opposite direction is not so easy: extending backwards from $x_k,x_{k+1},\ldots$ involves taking differences, which can produce values $\leq 0$.
We can show that entrance solutions exist by using a compactness argument
on the sequence of ratios between successive points:

\begin{theorem}
If $\pi$ satisfies (\ref{perm}), then there exists at least one solution of (\ref{pi}).
\end{theorem}

\proof
Given ${\bf x} \equiv (x_i)_{i \in \BZ} \in (0,\infty)^{\BZ}$,
%let $q_i := x_{i}/x_{i+1}$
%and let $V_i := (q_{i-1},q_{i-2},\ldots) $ $\in (0,\infty)^{\BN}$.
let $V_i := (x_{i-1}/x_i,x_{i-2}/x_{i-1},x_{i-3}/x_{i-2},
\ldots) $ $\in (0,\infty)^{\BN}$.
Note that if the sequence ${\bf x}$ satisfies
(\ref{pi}), then for all $i \in \BZ$,
$V_i \in [\frac{1}{2},1)^{\BN} \subset K := [\frac{1}{2},1]^{\BN}$.

If $V_i=(r_1,r_2,\ldots)$ then, factoring out $c := x_i$, we have
$(x_i,x_{i-1},x_{i-2},\ldots)$ $=c(1,r_1,r_1r_2,\ldots)$.  If
%$\pi(i)=i-d, \ d\geq 0$ so that
$x_{i+1}=x_i+x_{i-d}$ then
$x_{i+1}=c(1+\prod_{1\leq j \leq d} r_j) =: c/r_0$; this defines
$r_0$, such that  $V_{i+1} =(r_0,r_1,r_2,\ldots)$.
[In case $d=0$ the product is empty and has value 1, and
$x_{i+1}=2x_i, r_0=\frac{1}{2}$.]

This
motivates us to define, for each $d \geq 0$,
 a function $f^{(d)}\!\!: \ (0,\infty)^{\BN} \rightarrow  (0,\infty)^{\BN}$
by $(r_1,r_2,r_3\ldots) \mapsto ((1+\prod_{1 \leq j \leq d} r_j)^{-1},r_1,r_2,\ldots)$.
Note that these functions are continuous and map $K$ into itself.
Given $\pi$ satisfying (\ref{perm}), define $d(i)~:~=~i~-~\pi(i)$, and for
$n=1,2,\ldots$ define a map $T_{-n,0} := f^{(d(-1))} \circ f^{(d(-2))} \cdots \circ
f^{(d(1-n))}~\circ~f^{(d(-n))}$.  Note that if ${\bf x}$ satisfies (\ref{pi}) then
$T_{-n,0} (V_{-n}) = V_0$.
Conversely, given {\em any }
$\ldots, x_{-n-2},x_{-n-1},x_{-n} \in (0,\infty)$, (not necessarily satisfying (\ref{pi}),)
the map $T_{-n,0}$ gives
a recipe for extending the sequence with values $x_{-(n-1)},\ldots,x_{-1},x_0$
such that $x_{i+1}=x_i+x_{\pi(i)} \in (0,\infty)$ for $i=-n,\ldots,-2,-1$.
For $n=1,2,\ldots$ let $S_n$ be the image of the compact set $K$ under $T_{-n,0}$,
so $\cdots \subset S_2 \subset S_1 \subset K$. By the finite intersection
property, $\cap_{n\geq 1} S_n \neq \emptyset$.  Any point $(r_1,r_2,\ldots) \in \cap_{n\geq 1} S_n$ yields a sequence $x_0 =1, x_{-1}=r_1,$
$x_{-2}=r_1r_2,\ldots$ satisfying
$x_{i+1}=x_i+x_{\pi(i)}$ for $i=-1,-2,\ldots$, and this entrance solution
can be extended to a full solution of (\ref{pi}). \qed

Clearly if a sequence $(x_i)$ satisfies (\ref{pi}) then so
does any scalar multiple $(c x_i)$ for any $c>0$.  Are solutions
uniquely determined by $\pi$, up to such scalar multiples?
At least from the point of view taken in the proof of the
following theorem, the situation resembles time inhomogeneous
renewal chains; and we can only handle the bounded case.  Do
transient chains somehow correspond to counterexamples to
uniqueness?

\begin{theorem}
Assume $\pi$ satisfies (\ref{perm}) and $k := \sup \{ i - \pi(i): \ i \leq 0 \} < \infty$.
Then solutions of (\ref{pi}) are unique up to a scalar multiple.
\end{theorem}

\proof
We consider the effect of using $x_{i+1}=x_i+x_{\pi(i)}$ for $i=0,-1,\ldots,-n$
to express $x_1$ and $x_0$ as linear combinations, with nonnegative integer weights,
of $x_{-n},x_{-n-1},\ldots,x_{-n-k}$, viewed as indeterminates.
The idea is to show
that for large $n$, the weights in the combination for $x_1$ are close to being
a multiple of the weights for $x_0$, so that the ratio $x_1/x_0$ is close to being
determined by $\pi$.

The Hilbert
projective metric $\rho$ on $(0,\infty)^{k+1}$ (modulo scalar multiples)
is defined by $\rho(u,v) :=
\max_{i} \log (u_i/v_i) -\min_{i} \log (u_i/v_i)$.
If $u=(x_{-n},x_{-n-1},\ldots,x_{-n-k})$ comes from a solution of (\ref{pi}) and
$v=(1,1,\ldots,1)$, then the property $x_i < x_{i+1} \leq 2x_i$
implies that $\rho(u,v) \leq k \log 2$.

Consider $k+1$ by $k+1$ matrices with nonnegative integer coefficients,
indexed by $0 \leq i,j \leq k$, as follows.  Write $E^{(i,j)}$ for the matrix
having all zeroes, except for a single one in row $i$, column $j$.  Write
$B=\sum_{0<i\leq k} E^{(i,i-1)}$ for the matrix with ones below the
diagonal, and let $C=B+E^{(0,0)}$. For $0 \leq d \leq k$ let $A^{(d)} = C+ E^{(0,d)}$.
For $n \in \BZ$ let $d(n) =n- \pi(n)$,
$F^{(n)} = A^{(d(n))}$, and $M^{(n)} = F^{(-1)}  F^{(-2)}
 \cdots F^{(-n)}$.
Row $i$ of $M^{(n)}$ gives the coefficients of $x_{-i}$ as a linear combination
of $x_{-n},x_{-n-1},\ldots,x_{-n-k}$, and row 0 plus row $d(0)$ gives the
coefficients of $x_1$.

Note that  $C^k$ is all ones in column zero, and
zeroes elsewhere.  Any product $M$ of $k$ or more factors of the form
$A^{(d(i))}$ has all entries in column
zero strictly positive, and any product with
$2k$ or more factors has the property that every column is either all
zeroes, or else has all entries strictly positive. The maximum entry
in a product $M$ with exactly $2k$ factors is at most $s=2^{2k}$,
and the least entry in a non-zero column is at least one; hence if
$r=\rho(u,v)$, then
$$
 \rho(Mu,Mv)  \leq \log \left( \frac{1+se^r}{1+s} \right).
$$
It follows that using $r=k \log 2$ and given
$\epsilon>0$, we can pick $n$ so that for all $u,v$ with
$\rho(u,v) \leq r$ we have $\rho(M^{(n)}u,M^{(n)}v) < \epsilon$.
Using $v=(1,1,\ldots,1)$ and
the remarks at the ends of the previous two paragraphs,
the ratio $a :=(M^{(n)} v)_{d(0)} /  (M^{(n)} v)_0) $ satisfies
$|\log a - \log ((x_1-x_0)/x_0)| <\epsilon$, for any solution of
(\ref{pi}).  This show that any two solutions of (\ref{pi}) have
exactly the same ratio $x_0/x_1$.  The same argument shows that
the permutation $\pi$ determines the ratio $x_i/x_{i+1}$ for each $i \in \BZ$.
\qed

Problem \ref{prob 2}
 also includes solutions  placing probability $1/m$ on
each point in an  orbit of period $m>1$ of the
spacing transformation.  The simple deterministic solutions discussed above
are the case $m=1$, fixed points. For the
general case, there are $m$ distinct
 deterministic sequences; the random process picks each of these with probability $1/m$ each; and the spacings of the
$k^{th}$ sequence are all distinct and
as a set give the points of the $(k+1)^{st}$ sequence  mod $m$.

For example, for $m>1$ there are solutions of the form
$\p(X_i=x_i)=1/m$, where the deterministic sequence is
given by $x_i := b^i$, such that
$\Delta x_i := x_{i+1}-x_i$ $=(b-1)b^i$,
$\ldots, \Delta^m x_i = (b-1)^m x_i$ give $m$
distinct sets of points, and the last of these sets is
the same as $\{ x^i: \ i \in \BZ \}$.  This is possible iff
$(b-1)^m=b^{k}$ for some $k \in \BZ, b \in (1,\infty)$.
For example, with $m=2$ and $k=1$, $b=(3+\sqrt{5})/2 \doteq 2.618$
and with $m=2$ and $k=-1$, $b \doteq 1.75487$.
These examples with geometric sequences are misleading
in that in general, the $y_i := x_{i+1}-x_i$ are
not in increasing order, and the iterated transformation
on sequences does {\em not} correspond to $\Delta^m$, but rather
to (RANK $\circ \Delta)^m$.

To compare simple deterministic solutions with the scale invariant
Poisson processes, write  $\sigma$ for the inverse of $\pi$.  Note that for all $i$,
$i \leq \sigma(i)$ and $y_{\sigma(i)}=x_{i}<x_{i+1}=y_{\sigma(i+1)}$.
For the random solutions (\ref{random spacing})
there is a random permutation $\sigma$ of $\BZ$ such that for all $i$
$$
Y_{\sigma(i)} < Y_{\sigma(i+1)},
$$
with  $\sigma$  determined only up to translation.  If we write
$\sigma(i)=i+C(i)$, then it is easy to show, for the scale invariant Poisson process,
that a.s. $\limsup C_i=\infty$ and $\liminf C_i = -\infty$, a qualitative property
possessed by  no mixture of simple deterministic solutions.  Large deviations
for the permutation $\sigma$ are studied in \cite{zeitouni}.

To summarize the above discussion:  the solution to problem \ref{prob 2}
includes  \ a)  the scale invariant Poisson processes,
   b) simple deterministic solutions, i.e. fixed points
of the spacing transformation, \ c) deterministic orbits of
length $m>1$, and \ d) other extreme points of the set of distributional
solutions.  For b), it remains to resolve the question
of uniqueness relative to permutations satisfying (\ref{perm}), but with
unbounded displacements; and for c) and d), everything is open.

Peter Baxendale recently asked, for simple
deterministic solutions, what are the possible values of the ratio
$r:=x_0/x_1$ of two adjacent points?  Is $\frac{2}{3}$ achievable?
More generally,
writing $r_i := x_{i-1}/x_i$, what are the possible configurations
of $k+1$ consecutive points $(x_0,x_1,\ldots,x_k)$, i.e.~
what points in $[\frac{1}{2},1)^k$ are
realizable  as the value of $(r_1,r_2\ldots,r_k)$, for $k=1,2,\ldots$ \ ?
Two other ways to generalize the question about $r$ are to ask
which finite sets $A \subset (0,\infty)$ can satisfy
$A \subset \{ x_i\!: i \in \BZ \}$ for some  solution
of (\ref{simple det}), and similarly
which finite sets $B \subset [\frac{1}{2},1)$
can satisfy $B \subset \{ r_i\!: i \in \BZ \}$~?

\section{Problem session: a phase transition at $\ \theta=1 / \log 2$}

For any $\theta>0$, starting from a realization $\{ X_i: \ i \in \BZ \}$ of
the scale invariant Poisson process with intensity $\theta / x \ dx$, let
$A \equiv A(\theta)$ be the random closed set which is the closure of the countable
set whose points are $\sum_{i \in I} X_i$ for finite $I \subset \BZ$.
From the scale
invariance of the underlying $\{ X_i \}$ it follows easily that $A$ is also scale invariant:
for any $c>0$ and for each $\theta>0$
$$
c A(\theta) \indist A(\theta).
$$

The process $(A(\theta))_{\theta>0}$ has stationary, independent
increments in the following sense.
The Minkowski sum of two sets is $A \oplus B := \{ a+b: \ a \in A, b \in B \}$.
With $\X(\theta)$ to denote the set of points in the scale invariant Poisson process with intensity $\theta / x \ dx$,  the process $(\X(\theta))_{\theta>0}$ and hence also
$(A(\theta))_{\theta>0}$ can be constructed
with increasing sample paths:
$\theta_1<\theta_2$ implies $\X(\theta_1) \subset \X(\theta_2)$, and
hence $A(\theta_1) \subset A(\theta_2)$.
The process $\X$ has stationary, independent increments in the strong sense that $\X(\theta_2) \setminus
\X(\theta_1)$ is independent of $\X(\theta_1)$ and equal
in distribution to $\X(\theta_2-\theta_1)$, but the process
$A(\theta)_{\theta>0}$ has
stationary independent increments  in a weaker sense.
If $\theta_1<\theta_2$, and if $A'(\theta_2-\theta_1)$ is independent of
$A(\theta_1)$ and equal in distribution to $A(\theta_2-\theta_1)$, then
$A(\theta_1) \oplus
A'(\theta_2-\theta_1) \indist A(\theta_2)$.

Some  of the structure of the set of divisors of
a random integer, as in \cite{divisors}, is captured by the scale invariant closed set $A \subset [0,\infty)$.  First, for any $\theta > 0$ define a random
closed set $B(\theta) \subset [0,1]$ by $B := \{ \sum_{i \in I} V_i: \ {I \subset \BN} \}$.
From the Poisson-Dirichlet convergence in (\ref{billingsley}) and its
extension to the large deviation case for $\theta \neq 1$, it follows
easily that the random finite set
\begin{equation}\label{def D}
{\D}_n := \{ \log d /  \log n \ \} \subset [0,1],
\end{equation}
where $d$ runs over
the divisors of our random integer, has ${\D}_n \todist B(\theta)$.
(See Theorem \ref{divisors thm} for an extension.)  In
fact,
with the Hausdorff metric for closed subsets of $[0,1]$, and the $l_1$ metric
on $\BR^\infty$, the map which produces ${\D}_n$ from $(\log P_1,\log P_2,\ldots)$
and $B(\theta)$ from $(V_1,V_2,\ldots)$ is a contraction.  Second, just like
(\ref{blow up}) it is easy to see that as $v \ra \infty$, $v B(\theta)  \todist A(\theta)$.
In fact, for any $0<x \leq v$,
\begin{equation}\label{blow up B}
d_{TV}(v B \cap (0,x), A \cap (0,x) \ )
 = d_{TV}( B \cap (0,x/v), (v^{-1}A) \cap (0,x/v) \ )
\end{equation}
$$
\leq  d_{TV}( \V \cap (0,x/v), (v^{-1}\X) \cap (0,x/v) \ )
=H_\theta(x/v).
$$
Thus for fixed $x$, $d_{TV}(v B \cap (0,x), A \cap (0,x) \ )  \ra 0$
as  $v \ra \infty$, which is even stronger than the distributional convergence
$v B \todist A$.

Let $f(\theta) := \p(1 \in A(\theta)\ )$. Using scale invariance,
$\forall x>0, \p(x \in A(\theta)\ ) = f(\theta)$. Hence with $m(A)$ to
denote Lebesgue measure,  for $0 \leq a < b < \infty$,
\mbox{$\e m(A~\cap (a,b) \, )$} $ = (b-a) f(\theta)$.  It is easy to show that for $\theta  \leq 1/ \log 2$,
$f(\theta)=0$.  Note that for any $\theta$, $f(\theta) \leq \p(T<~1)$ $ < 1$.
In lecture, we asked a semi-open question:
\begin{problem} \label{problem log 2}
Prove the conjecture that if
$\theta > 1 / \log 2$, then $f(\theta) > 0$.
\end{problem}
At the time of lecture, the conjecture was semi-open in the
following sense: Tenenbaum (private communication) had proved a
related property about the set of divisors of a random integer from
1 to $n$ with $\theta \log \log n $ distinct prime divisors.
From the early version of this number theory result,
it was not yet possible to deduce the continuum result as a corollary, so the conjecture
was indeed open.
However,  it was highly plausible to guess that
the underlying principles from Tenenbaum's proof would also work directly
on the continuum problem.  Subsequent to the workshop, the conjecture has
indeed been proved this way.
Further, Tenenbaum gave a sharper version of the number theoretic result,
from which the continuum result  follows as a corollary.  This amounts
to using a more complicated process, prime divisors of a random integer,
to approximate a simpler process, the scale invariant Poisson!

The event that $A(\theta)$ has positive length is a tail event with respect to the
independent exponentials $W_i$ in (\ref{u from w}), and hence by the
Kolmogorov zero-one law and the above, for every $x>0$ and
any $\theta > 1/ \, \log 2$, $1=\p( m(A \cap (0,x) \, ) > 0 )$.

The random sets above are closely related to the theory of Bernoulli
convolutions; we learned of the connection thanks to Jim Pitman.
For Bernoulli convolutions, the usual setup is to start with a deterministic
sequence $r_1,r_2,\ldots > 0$ with $\sum r_n < \infty$
 and to consider
 the random variable $Y = \sum_{i \geq 1} S_i r_i$, where the $S_i$ have values 1, -1, independently with probability 1/2 each.
The closed support of the distribution of $Y$ is the set $K$
of all points of the form $\sum s_i r_i$, where the $s_i$ are arbitrarily chosen from
$\{1, -1\}$.  There is a straightfoward
translation to the situation where $s_i \in \{0,1\}$, which is natural
for the application to number theory.

Erd\H os focussed on the case $r_n := \lambda^n$ for a fixed $0<\lambda<1$.
Here it is obvious that $K$ has zero length for $0<\lambda<\frac{1}{2}$,
and that $K$ is an interval for $\frac{1}{2} \leq \lambda < 1$.
Erd\H os asked, when is the distribution $\mu_\lambda$
of $Y$ absolutely continuous with respect to Lebesgue measure?
For $\lambda<\frac{1}{2}$,
since the support $K$ has zero length, $\mu_\lambda$ is singular, but for $\lambda>\frac{1}{2}$,
even though the support $K$ is an interval, the absolute continuity question is subtle.  In  1939 Erd\H os \cite{erdos39}
showed that there are values in  $(\frac{1}{2},1)$, such as $\lambda =  2/(1+\sqrt{5})$,
for which
$\mu_\lambda$ is singular, and in 1940 \cite{erdos40} he
showed that there exists some $t<1$ such that for almost every $\lambda \in (t,1)$,
$\mu_\lambda$ is absolutely continuous.  In 1995,  Solomyak \cite{solomyak}
 showed that $t$ in the previous sentence can be taken to be 1/2.  Peres and Solomyak
\cite{peres and solomyak}  give a simple proof of this.
In 1958 Kahane and Salem \cite{kahane58}
discussed a variant of the problem where the summands $r_n$ are random,
of the form $r_n = U_1 U_2 \cdots U_n$ where the $U_i$ are independent,
and uniformly chosen from intervals $[a_i,b_i] \subset [\frac{1}{2},1]$.

Consider the random closed
set $C := \{ \sum_{i \geq 1} c_i X_i: \ c_i \in \{0,1\} \ \}$,
with the scale invariant Poisson process labelled as in (\ref{at zero}), so that $C \subset [0,T]$.  This is close to $A$ in that always $A \cap [0,1) = C \cap [0,1)$, and related to
$B$: as a corollary of (\ref{t1}), for each $\theta > 0$, $B \indist (C \ | \ T=1 \ )$.
Let  $J_1,J_2,\ldots$ be fair coins with values in $\{0,1\}$, independent of
each other and everything else.  The random series  $Y^ *:= \sum_{i \geq 1} J_i X_i$
has values in the random set $C$, and the random series
$Y := \sum_{i \geq 1} J_i V_i$ has values in the random set $B$.  Consider the conditional
distribution $\mu_{Y^* | X}$ of $Y^*$ given $X_1,X_2,\ldots$, and the conditional distribution $\mu_{Y|V}$ of $Y$ given the Poisson-Dirichlet process.  These
are {\em random} probability measures on $[0,\infty)$ and $[0,1]$,
respectively, and  correspond to
the fixed probability  distribution $\mu_\lambda$ from the
situation where only the signs are random.
The method of \cite{peres and solomyak} can be adapted to show that
for $\theta > 1/\log 2$, a.s. the distributions $\mu_{Y^*|X}$ and
$\mu_{Y|V}$ are absolutely continuous with respect to Lebesgue
measure, with density in $L_2$ (Yuval Peres, private communication).

The closed support in $[0,\infty)$ of $\mu_{Y^*|X}$ is $C$, and that of
$\mu_{Y | V}$ is $B$.  The distribution of $Y$ itself, without conditioning
on the value of the Poisson-Dirichlet process, is the probability
measure $\e \, \mu_{Y|V}$, averaging over the
values of $(V_1,V_2,\ldots)$.  From Hirth \cite{hirth} and Donnelly-Tavar\'e
 \cite{dt}, the distribution of $Y$
 is simply Beta($\frac{\theta}{2},\frac{\theta}{2}$), which
for $\theta=1$ is the arc-sine distribution.  It gives the limit of the
distribution of $\log d/ \log n$, and the limit distribution of
$\log d / \log N$, where a random integer $N$ is chosen uniformly from
1 to $n$, and then $d$ is chosen uniformly from the  $\tau(N)$ divisors of $N$;
see \cite{tenenbaum},  II.6.2.
All this background suggests looking at the distribution of (the rescaled log of )
a randomly chosen divisor of uniformly chosen random integer, as a random
probability distribution.  We need some notation.

Let $\P$ denote the space of probability measures on $\BR$, with the
topology of weak convegence, as given by the L\'evy metric $\rho$; see
e.g.~\cite{dudley}.
For $m \geq 1$ let $\mu_m := \left( \sum_{d|m} \delta_{\log d /  \log  m} \right) /
\left( \sum_{d|m} 1 \right) \ \in \P$ be the probability measure on [0,1]
that puts mass $1/\tau(m)$ on each point of the form $\log d / \log m$,
where $d$ runs over the divisors of $m$.  (For $m=1$, interpret 0/0 as 1,
so that $\mu_1 := \delta_1$.)

For $n \geq 1$ let $N \equiv N(n)$ be chosen uniformly from
1 to $n$, i.e.~ $\p_n(N=m)$ $ =1/n$ for $m=1,2,\ldots, n$.
Let $\mu_N$ denote the random measure whose value is
$\mu_m$ on the event $\{ N =m \}$, so that the law
of $\mu_N$ as a random element of $\P$ is
$(1/n) \sum_{1 \leq m \leq n} \delta_{\mu_m}$.

\begin{theorem}\label{divisors thm}
As $n \ra \infty$,
$$
\mu_N \todist \mu_{Y|V},
$$
i.e.~as random a element of $\P$,  $\mu_{N(n)}$
converges in distribution to $\mu_{Y|V}$, using the Poisson-Dirichlet
process with $\theta=1$.
\end{theorem}

\proof
As noted in \cite{hirth}, if $m = p_1 p_2 \cdots p_k$ is
squarefree, then $\mu_m$, supported at $2^k$ distinct points,
is a Bernoulli convolution, and the essential
difficulty is to deal with $m$ such that $p^2|m$ for some
prime.  Write $P_i(m)$ for the  $i^{th}$ largest prime factor of
$m$, with $P_i=1$ if $i$ is greater than
the number $\Omega(m)$ of prime factors
of $m$, including multiplicities.
Let $J_i$ be $\{0,1\}$-valued fair coins, independent of
everything else.
For $k=1,2,\ldots$ define  $\mu_m^{(k)} \in \P$ to be the
distribution of $\sum_1^k J_i \, (\log P_i(m) / \log m)$, so
that $\mu_m^{(k)}=\mu_m$ iff $m$ is squarefree and has
at most $k$ prime factors. Similarly, define $\mu_{Y|V}^{(k)} \in \P$
to be the conditional
distribution of $\sum_1^k J_i V_i$ given $(V_1,V_2,\ldots)$.
Write $\L(S)$ for the law of $S$.
The L\'evy metric $\rho$ on
$\P$ has the property that for random variables $S,T, \,$ $1=\p(|S-T|\leq \epsilon)$
implies that
$\rho(\L(S),\L(T)) \leq\epsilon$.
Hence with  $S=\sum_1^k J_i v_i$ and
$T=\sum_1^\infty J_i v_i$, it follows
 that $\rho(\mu_{Y|V}^{(k)},\mu_{Y|V})
$ $ \leq V_{k+1}+V_{k+2}+\cdots$ $ = 1-(V_1+\cdots + V_k)$.

Given $\epsilon>0$, pick $k$ so that $\p(V_1+\cdots+V_k < 1-\epsilon)
< \epsilon$.  This
implies $\p(\, \rho(\mu_{Y|V}^{(k)},\mu_{Y|V})>\epsilon)<\epsilon$.
For any $m$ write $r= P_1(m)\cdots P_k(m)$ and   consider the map $g_m$ with $g_m(d) = (d,r)$, the greatest common divisor.  This map, applied
to the set of divisors of $m$, is exactly $\tau(m)2^{-k}$ to one,
{\em provided that} $m$ has
at least $k$ prime divisors and the $k$ largest prime divisors are distinct.  For
$m$ satisfying this condition, $g_k$ applied to a divisor chosen randomly
from the $\tau(m)$ divisors of $m$ yields a uniform pick from
the $2^k$ divisors of $r$, and hence $\rho( \mu_m^{(k)},\mu_m)
\leq \log(m/r)/\log m$.

We note that Billingsley \cite{billing large} actually proved both (\ref{billingsley})
and the slightly different form we need here, that $(\log P_i(N)/\log N)_{i \geq 1}
$ $\todist (V_i)_{i\geq 1}$.  Thus we can pick $n_1$ so that for $n \geq n_1$,
$\p_n(P_1(N)\cdots P_k(N) \geq N^{1-\epsilon}$ and $  P_1(N)> \cdots >P_k(N)>1)$ $> 1-\epsilon$.
This yields, for $n \geq n_1$, that $\p_n(\, \rho(\mu_N^{(k)},\mu_N) > \epsilon) < \epsilon.$

The map $f_k\!: (x_1,\ldots,x_k) \mapsto \L( \sum_1^k J_i x_i)$,
from $(\BR^k,l_1)$ to $(\P,\rho)$  is continuous;
in fact it is a contraction.
Now  $\mu_{Y|V}^{(k)} = f_k(V_1,\ldots,V_k)$, and
$\mu_N^{(k)}=$ \mbox{$f_k(\log_N P_1(N),\ldots,\log_N P_k(N))$}, so
using Billingsley's result again,
as $n \ra \infty$, $\mu_{N}^{(k)} \todist \mu_{Y|V}^{(k)}$.  By
 Skorohod embedding \cite{dudley} for the complete separable metric space $\P$,
there exists a coupling in which $ \mu_{N(n)}^{(k)} \ra \mu_{Y|V}^{(k)}$
almost surely.  This coupling can be extended to a coupling of $N(1),N(2),\ldots$
and the Poisson-Dirichlet process, so that for one fixed $k$,
$\mu_{N(n)}^{(k)}, \mu_N(n),
 \mu_{Y|V}^{(k)}$, and $\mu_{Y|V}  $  are all
defined on the same probability space,
and a.s. $ \mu_{N(n)}^{(k)} \ra \mu_{Y|V}^{(k)}$.  Pick
 $n_2$ such that for $n \geq n_2$,
$\p (\, \rho(\mu_{N(n)}^{(k)}, \mu_{Y|V}^{(k)})>\epsilon)<\epsilon$.
The triangle inequality yields, for $n\geq \max(n_1,n_2)$, that
$\p(\, \rho(\mu_N(n),\mu_{Y|V}) > 3 \epsilon) < 3 \epsilon$. \qed

%\section{References}

\end{document}